\documentclass[11pt]{amsart}
\usepackage[english]{babel}
\usepackage{graphicx}
\usepackage[T1]{fontenc}
\usepackage{color}
\usepackage{bigints}
\usepackage{soul}
\usepackage{fullpage}
\newtheorem{Ass}{Assumption}
\newtheorem{Thm}{Theorem}
\newtheorem{Prop}[Thm]{Proposition}

\newtheorem*{Prop*}{Proposition}
\newtheorem*{Cor*}{Corollary}
\newtheorem*{Thm*}{Theorem}
\usepackage{hyperref}
\usepackage[utf8]{inputenc}
\usepackage[english]{babel}
\usepackage{etoolbox}
\apptocmd{\sloppy}{\hbadness 10000\relax}{}{}
\hypersetup{
	colorlinks=true,
	linkcolor=blue,
	urlcolor=cyan,
}

\begin{document}
	\title[Existence of solution of a triangular degenerate reaction-diffusion system]{Existence of solution of a triangular degenerate reaction-diffusion system}

	\author[Saumyajit Das]{Saumyajit Das}
	\address{Department of Mathematics, Indian Institute of Technology Bombay, Powai, Mumbai 400076 India.}
	\email{194099001@iitb.ac.in}
	\begin{abstract}
		In this article we study a chemical reaction-diffusion system with $m$ unknown concentration. The non-linearity in our study comes from a particular chemical reaction where one unit of a particular species generated from other $m-1$ species and disintegrates to generate all those $m-1$ species in the same manner, i.e., triangular in nature. Our objective is to find whether global in time solution exists for this system where one or more species stops diffusing. In particular weak global in time solution exists for all the degenerate cases in any dimension. We are further able to show classical global in time solution exists for all the degenerate cases in any dimension except one and this particular case too attain classical global in time solution upto dimension $2$. We also analyze global in time existence result for the case of quadratic non-linear rate functions and also analyze a three dimensional case.
	\end{abstract}
	
	\maketitle
	
	\section{Introduction}\label{sec:introduction}
		Reaction-Diffusion system, originating from chemical kinetics,governs the evolution  of concentration of  species in time at various spatial locations. The nature of such system is the species are simultaneously diffusing and undergoing chemical reactions i.e., reversible in nature. 
		
		
		
		
		In this article we consider a $m$ species triangular reaction-diffusion system where in a domain $m$ species  $ X_1,\cdots,X_{m}$ are undergoing through the following reversible system:
		
		\[
		\alpha_1 X_1+\alpha_2X_2+\cdots+\alpha_{m-1}X_{m-1} {\rightleftharpoons} X_{m}
		\]

		The spatial domain is taken to be a  bounded domain $\Omega\subset\mathbb{R}^N$ with $C^{2+\nu}$ boundary. For the unknown $a_i:[0,T)\times\Omega\to\mathbb{R}$, represents the concentrations species $X_i$, for all $i=1,\cdots,m$. The reaction-diffusion system is given by the following system of partial differential equations
		
		\begin{equation}\label{triangular  model}
			\left \{
			\begin{aligned}
				\partial_t a_i - d_i \Delta a_i=& a_m-\displaystyle\prod_{1}^{m-1}(a_j)^{\alpha_i} \qquad \qquad \mbox{in}\ (0,T)\times\Omega\\
				\partial_t a_m-d_{m}\Delta a_m=& \displaystyle\prod_{1}^{m-1}(a_j)^{\alpha_i}-a_m \qquad \qquad \mbox{in}\ (0,T)\times\Omega \\
				n.\nabla_{x}a_i= & 0 \qquad \qquad \qquad \qquad \qquad \ \ \mbox{on}\ (0,T)\times \partial \Omega \\
				a_{i}(0,x)=& a_{i,0}(\in C^{\infty}(\Omega))\geq 0,  
			\end{aligned}
			\right .
		\end{equation}
		
	where $n(x)$ denotes the outward unit normal to $\Omega$ at the point $x\in\partial\Omega$. Furthermore  $\alpha_i$ is the stoichiometric coefficient corresponds to the species $a_i$, generally taken as non-negative integer and the diffusion coefficients $d_i$ are taken to be non-negative real number.

	\vspace{.2cm}
	
	A reaction-diffusion system is called triangular system  if there exists some positive lower triangular matrix whose action with the source terms makes it linearly dominated by the unknowns. The rate functions corresponding to our system
	
	\begin{equation*}
		\left \{
		\begin{aligned}
		f_i = & a_m-\displaystyle\prod_{1}^{m-1}(a_j)^{\alpha_i} \quad \forall i=1,\cdots,m-1\\
		f_m=& -f_1.
		\end{aligned}
		\right . 
	\end{equation*}

     Consider the positive lower triangular matrix
     \begin{equation*}
     	P=(p_{ij})_{m\times m}= \left \{
     	\begin{aligned}
     		&1 \qquad i=j, i\in\{1,\cdots,m \}\\
     		&1 \qquad j=i-1, i\in \{2,\cdots,m\}\\
     		& 0 \qquad \text{otherwise}.
     	\end{aligned}
     	\right .
     \end{equation*}
 
     The action of this matrix on the rate function vector $(f_1,\cdots,f_m)$ is the the following product
     
     \[
     P \begin{pmatrix}
     	f_1 \\ \cdot \\ \cdot\\ \cdot \\ f_m
     \end{pmatrix} \leq \left[1+\sum\limits_{1}^{m}a_i\right] \begin{pmatrix}
     1 \\ 2\\ \cdot \\ \cdot \\ 2 \\ 0
 \end{pmatrix},
     \]
	 which makes our system a triangular one.
	    
	    When all the diffusion coefficients are strictly positive i.e., $d_i>0$ for all $i=1,\cdots,m$(here onwards we call this as non-degenerate setting),
	    well posedness of this triangular system and global in time existence of classical solution  is established in  articles \cite{Pie2010} \cite{Pierre2003}. In the article \cite{DF15}, authors first analyze the simplest three species triangular reaction-diffusion system where one or more species stops diffusing(here onwards we call this as degenerate setting). Authors consider the following reversible reaction:
	    \[
	    X_1+X_2 \rightleftharpoons X_3,
	    \] 
	    
	    where the concentrations $a_1,a_2,a_3$ satisfy the following system 
	    
	    \begin{equation*}
	    	\left \{
	    	\begin{aligned}
	    		\partial_t a_1 -d_1 \Delta a_1 =& a_3-a_1a_2 \qquad \mbox{in}\ (0,T)\times\Omega\\
	    		\partial_t a_2 -d_2 \Delta a_2 =& a_3-a_1a_2 \qquad \mbox{in}\ (0,T)\times\Omega\\
	    		\partial_t a_3 -d_3 \Delta a_ =& a_1a_2-a_3 \qquad \mbox{in}\ (0,T)\times\Omega\\
	    		n.\nabla_x a_i=& 0 \qquad \qquad \mbox{on}\ (0,T)\times\partial\Omega, \forall i=1,2,3\\
	    		a_{i}(0,x)=& a_{i,0}(\in C^{\infty}(\Omega))\geq 0 \qquad \forall i=1,2,3.
	    	\end{aligned}
	    	\right .
	    \end{equation*}

	      There are following kinds of degeneracies:
	      
	      \begin{itemize}
	      	\item $d_1=0$, $d_2$,$d_3>0$ or $d_2=0$,$d_1$,$d_3>0$,
	      	\item $d_3=0$,$d_1$,$d_2>0$,
	      	\item $d_1=d_2=0$, $d_3>0$,
	      	\item $d_1=d_3=0$, $d_2>0$ or $d_2=d_3=0$, $d_1>0$.
	      \end{itemize} 
      
      In the article  \cite{DF15}, it has been shown except for the case $d_3=0,d_1,d_2>0$, system always has smooth global in time solution. However and for the degenerate case $d_3=0,d_1,d_2>0$ one can find global classical solution upto 3 dimension, although global in time weak  $L^1((0,T)\times\Omega)$ solution exists in any dimension. It requires various duality arguments as described in \cite{CDF14} \cite{DFMV07}\cite{Amann1985}\cite{rothe06}. The following two results one of which ensures positivity of solution and another lays the idea of maximal time interval, where we can find classical solution, are extremely useful.
      
      	\begin{Thm*}[Positivity of solution] 
      		 Let $ u_i:(0,T)\times\Omega\rightarrow \mathbb{R}$, satisfies the following equation(weak sense too included) with $d_i\geq 0$
      	
      	\begin{equation*}
      		\left \{
      		\begin{aligned}
      			\partial_t u_i- d_i \Delta u_i= & f_i(u_1,\cdots,u_i) \qquad \qquad \qquad \mbox{in}\ (0,T)\times\Omega, \forall i=1,\cdots,m \\\
      			n.\nabla_{x} u_i = & 0 \qquad \qquad  \qquad \qquad \qquad \ \ \  \ \mbox{on}\ (0,T)\times\partial \Omega\\
      			 u_i(0,x)& \in C^{\infty}(\Omega)\geq 0.
      		\end{aligned}
      		\right .
      	\end{equation*}
      	Furthermore let $f=(f_1,\cdots,f_m):\mathbb{R}^n\rightarrow \mathbb{R}^m$ is quasi-positive i.e.,
      	\[
      	f_i(r_1,r_2,\cdot \cdot,r_{i-1},0,r_{i+1},\cdots,r_n) \geq 0, \  \forall (r_1,r_2,\cdots,r_n)\in [0,+\infty)^n.
      	\]
      	 Then the solution remains non-negative i.e., $u_i \geq 0, \ \forall i=1,\cdots,m$. 
      	\end{Thm*}
      	 
      	 Consider triangular reaction-diffusion equation \eqref{triangular  model}. Rate functions corresponding to this system satisfy above quasi-positive condition. Hence if solution exists then it will be non-negative provided the initial condition is non-negative \cite{Pie2010}\cite{Pierre2003}\cite{Amann1985}.

      	\vspace{.2cm}
      	
      	\begin{Thm*}[Maximal interval of solution]
      	 Let $u :(0,T)\times\Omega\rightarrow \mathbb{R}$, satisfies the following equation(weak sense is also included) with $d\geq 0$ and $f\in C^1(\mathbb{R})$:
      	
      	\begin{equation*}
      		\left \{
      		\begin{aligned}
      			\partial_t u- d \Delta u= & f(u) \qquad \qquad \qquad \ \ (0,T)\times\Omega\\
      			n.\nabla_{x} u = & 0 \qquad \qquad \qquad \qquad  (0,T)\times\partial \Omega\\
      			u(0,x) \in C^{\infty}& (\Omega). 
      		\end{aligned}
      		\right .
      	\end{equation*}
      	
      	if $T_{max}$ denotes the maximal interval of existence then
      	
      	\[
      	\lim \limits_{t\uparrow T_{max}} \sup\limits_{\Omega}\vert u\vert < +\infty \ \text{and vice versa}.
      	\]
        
      	This indicates, to show global in time existence we need to show the space-time supremum norm bound of a species at any particular time is bounded i.e.,  showing $\lim \limits_{t\uparrow T} \sup\limits_{\Omega}\vert u\vert < +\infty$ is enough. Proof of this theorem can be found in the article \cite{Amann1985}.
      \end{Thm*}
   This idea is followed  widely to show global in time existence of solution \cite{Pie2010}\cite{DF15}\cite{FMT20}\cite{DFMV07}\cite{EMT20}. In particular global in time smooth solution for non-degenerate triangular system has shown in the article \cite{Pie2010}. Similar for three species degenerate model has shown in the article \cite{DF15}.

	
	\vspace{.2cm}

	\vspace{.3cm}
	
	In this article we will discuss global in time existence of positive solutions of various degenerate triangular reaction-diffusion equation. We will show for $d_m>0$ or for $d_m=0$ alongwith some other $d_i=0$ system-\eqref{triangular  model} always have global in time smooth solution with the help of various duality estimates and parabolic regularity. Furthermore for the degeneracy $d_m=0$ and all other diffusion coefficients strictly positive,  we are able to show global in time weak $(H^1((0,T)\times\Omega_T))^{m-1}\times L^1((0,T)\times\Omega)$ solution exists for any dimension. However for this particular degeneracy we will show global in time  classical solution exists upto dimension $2$. We further discuss two special cases one on quadratic bound on the rate function and another is a special 3 dimensional case. Many of our ideas are borrowed from \cite{Pie2010}\cite{DF15}\cite{rothe06}\cite{FMT20} \cite{EMT20} \cite{Tan18}.
	
	\vspace{.2cm}
	
	We divide all the triangular degenerate cases in three different classes:
	\begin{flalign}
		\text{System}-\text{A}_1:& \quad d_m>0,d_i=0;\ \text{for some}\ i\neq m. \label{A1}\\
		\text{System}-\text{A}_2:& \quad      d_m=d_i=0,\ \text{where}\ i\neq m. \label{A2}\\
		\text{System}-\text{A}_3:& \quad d_m=0, d_i>0; \forall i=1,\cdots,m-1. \label{A3}
	\end{flalign}
	
	Furthermore we divide the diffusion coefficients into two different index classes basically one index class corresponds to diffusing species i.e., where the diffusion coefficient is non-zero and other corresponds to the non-diffusing species i.e., where the diffusion coefficient is zero.  We denote it by the following notation
	\[
	\Lambda_1= \{i: d_i=0\}, \ \Lambda_2=\{i:d_i>0\}.
	\]
	
	\vspace{.5cm}
	
	\textbf{\large{Some Further Notation}}:
	\begin{itemize}
		\item $\Omega_{\tau}= (0,\tau)\times \Omega$, \quad $\forall 0<\tau\leq T$.\\
		\item $\alpha_m=1$.
	\end{itemize}
	
	\vspace{.3cm}
	
	For all these three different classes of degenerate systems we first introduce one common approximate system where our rate functions will be smooth,Lipschitz and  the approximate system carries the same quasi-positive character. Hence for the approximate system  global in time smooth positive  solution exists provided the initial condition is smooth and positive\cite{rothe06}\cite{taylor2013partial}. The approximate system is described below
	
	\begin{equation}\label{approximate  system}
		\left\{
		\begin{aligned}
			\partial_t a^n_i - d_i \Delta a^n_i=& \frac{a^n_m-\displaystyle\prod_{j=1}^{m-1}(a^n_j)^{\alpha_j}}{\phi^n}  \qquad \qquad \ \mbox{\ in \ }(0,T)\times\Omega, \forall i=1,\cdots,m-1 \\
			\partial_t a^n_m-d_m\Delta a_m^n=& \frac{\displaystyle\prod_{j=1}^{m-1}(a^n_j)^{\alpha_j}-a^n_m}{\phi^n} \qquad \qquad \ \mbox{\ in \ }(0,T)\times\Omega \\
			n.\nabla_{x}a^n_i= & 0 \qquad \qquad \qquad \qquad \qquad \ \  \mbox{ \ on \ } (0,T)\times \partial \Omega \\
			a^n_{i}(0,x)=& a_{i,0}(\in C^{\infty}(\Omega))\geq 0,  
		\end{aligned}
		\right.
	\end{equation}
	where 
	\[
	\phi^n:=1+\frac{1}{n}\Big(\ \sum_{i=1}^{m}a^{n}_{i} \Big)^{Q+2} \quad \forall  n\in\mathbb{N}, \ \  Q=1+\sum\limits_{i=1}^{m-1}\alpha_{i}.
	\]

	We obtain a solution of the degenerate triangular reaction-diffusion system \eqref{triangular  model} as a weak limit of the solutions of the approximate system. We use various duality estimates motivated from the article \cite{Pie2010}. The key steps are described below.
	
	\vspace*{.2cm}
	
	Step 1: We extract a weak solution for those species which are diffusing i.e., non-degenerate species, which is the weak limit of the approximate solutions thanks to property of heat kernel and Dunford-Pettis' theorem.
	
	\vspace*{.2cm}
	
	Step 2: We will  find some uniform norm bound of the solutions of a particular species of the approximate system in some large $L^p$ space where the bound will be uniform with respect to the index $n$.
	
	\vspace*{.2cm}
	
	Step 3: By duality we will find uniform norm bound of the solutions of all  species of the approximate system in same large $L^p$ space where the bound will be uniform with respect to the index $n$.
	
	\vspace*{.2cm}
	
	Step 4: Using integrability estimate [see appendix \eqref{eq:Linfty-estimate}] we will arrive at a uniform $L^{\infty}$ time-space bound for all the species where the bound is uniform with respect to the index $n$. In the light of "maximal interval of solution" theorem, this estimate will establish global in time existence of solution.
	
	\vspace*{.2cm}
	
	Step 5: The approximate system corresponds to a degenerate species satisfies a sequence of o.d.e.'s, instead of p.d.e.'s. As the unknown also presents in rate functions, by Picard type iteration and Vitali convergence lemma we can extract  a weak solution for the ordinary differential equation corresponding to degenerate species. To  obtain this weak limit from the sequence of o.d.e.'s, we need  some condition on stochiometric coefficients. The particular condition is stated below
	
	\begin{Ass}[Stochiometric Condition(SC)]\label{SC}
		\begin{equation*}
			\left \{
			\begin{aligned}
				\alpha_i\in & [1,\infty), \forall i\in\Lambda_1,\\
				\exists \ \alpha_j\in & \{1\}\cup[2,\infty), \text{\ for some \ }j\in\Lambda_1.
			\end{aligned}
			\right .
		\end{equation*}    
	\end{Ass}
	\textbf{Note:} For the degenerate class system-$\text{A}_3$\eqref{A3}
	the stochiometric condition\eqref{SC} is automatically satisfied. So we don't need this condition for any result related to system-$\text{A}_3$.
	
	\vspace{.3cm}

	 The above mentioned iterative argument can be found in the article \cite{DF15}, where the authors analyze the existence result of a particular four species degenerate reaction-diffusion system. 
	 
	 \vspace{.2cm}
	 
	 For system-$\text{A}_3$\eqref{A3} step 1 and step 5 directly yields a global in time positive weak solution. However the step 3, we are able to obtain it upto 2 dimension thanks to the heat kernel estimate by fundamental solution of the heat equation as described in the articles \cite{Morra83}\cite{ML15}.


	\vspace{.5cm}
	
	Next we state our main results.
	
	\begin{Thm}\label{Existence A1}
		Under condition SC\eqref{SC}, there exists unique positive global in time  classical solution in any dimension for system-$\text{A}_1$\eqref{A1}. 
	\end{Thm}
	
	\begin{Thm}\label{Existence A2}
		Under condition SC\eqref{SC}, there exists unique positive global in time  classical solution in any dimension for system-$\text{A}_2$\eqref{A2}.
	\end{Thm}
	
	\begin{Thm}\label{Existence A3}
		There exists  positive global in time  weak $(H^1(\Omega_T))^{m-1}\times L^1(\Omega_T)$ solution for system-$\text{A}_3$\eqref{A3}, in any dimension. However  unique positive global in time  classical solution exists for system-$\text{A}_3$\eqref{A3} in dimension 1 and 2.
	\end{Thm}

	\section {\textbf{\large{Existence of solution of degenerate triangular system-$\textbf{A}_1$}}} 
		Our first goal is to extract a weak convergent subsequence of non-degenerate species from the solution of approximate system.
	
	Rate functions associated with each of the approximate system\eqref{approximate  system} is globally Lipschitz with positive smooth initial condition. Hence we have positive smooth global in time solution  $a_i^n:\Omega_T\rightarrow \mathbb{R}$, for each of the approximate system. Positivity comes from quasi-positive nature of each of the approximate system \cite{Amann1985} \cite{rothe06}\cite{quittner2019superlinear}.

	Let's define a positive functional based on the solutions of the approximate system
	
	\[
	E(a_i^n:i=1,2\cdots,m)=\int_{\Omega} \sum_1^m \alpha_i(a_i^n(\ln{a_i^n}-1)+1).
	\]
	
	Differentiating entropy with respect to time and replacing time derivative of the species from the approximate system\eqref{approximate  system}, we obtain:
	\[
	\frac{\mbox{d}E(a_i^n:i=1,2\cdots,m)}{\mbox{d}\ t}= - \int_{\Omega}\sum\limits_{i=1}^{m} \alpha_i d_i \frac{\vert \nabla a_i^n \vert^2}{a_i^n}\nonumber  -\bigint_{\Omega}  \frac{ \Big(a^n_m-\displaystyle \prod\limits_{j=1}^{m-1} (a^n_j)^{\alpha_j}\Big)}{\phi^n}\ln\left({\frac{a^n_m}{\prod\limits_{j=1}^{m-1} (a^n_j)^{\alpha_j}}}\right).
	\]

	Integrating the above relation with respect to time yields
	
	\begin{flalign*}
		\sup\limits_{t\in[0,T]}\int_{\Omega}&  \sum_1^m  \alpha_i(a_i^n(\ln{a_i^n}-1)+1)+ \int_{\Omega_T}\sum\limits_{i=1}^{m} \alpha_i d_i \frac{\vert \nabla a_i^n \vert^2}{a_i^n}\nonumber \\ & +\bigint_{\Omega_T}  \frac{ \Big(a^n_m-\displaystyle \prod\limits_{j=1}^{m-1} (a^n_j)^{\alpha_j}\Big)}{\phi^n}\ln\left({\frac{a^n_m}{\prod\limits_{j=1}^{m-1} (a^n_j)^{\alpha_j}}}\right)\leq E(a_{i,0}:i=1,2\cdots,m). 
	\end{flalign*}
	
	The following two bounds follows from the above relation
	
	\begin{equation}\label{Uniform-integrability traingular}
		\left\{
		\begin{aligned}
			\int_{\Omega} \sum\limits_{i=1}^m a^n_i \leq & M_2, \\
			\int_{\Omega} \sum\limits_{i=1}^m  \vert a_i^n\ln{a_i^n}\vert  \leq & M_2+E(a_{i,0}:i=1,2\cdots,m) )+ e^{-1}\vert\Omega\vert>0, 
		\end{aligned}
		\right .
	\end{equation}
	
	where, $M_2=\max\limits_{i\in\{1,\cdots,m\}}\{\alpha_i\}e^2\vert\Omega\vert+ \frac{\max\limits_{i\in\{1,\cdots,m\}}\{\alpha_i\}}{\min\limits_{i\in\{1,\cdots,m\}}\{\alpha_i\}}E(a_{i,0}:i=1,2\cdots,m)>0$.
	
	Furthermore condition \eqref{Uniform-integrability traingular} implies  uniform integrability of  $\sum\limits_{i=1}^m a_i^n,  \ \ \forall n\in \mathbb{N}$.  

	Consider the following algebraic inequality for two positive quantity $x,y$
	\[
	x \leq \kappa y +\frac{1}{\ln \kappa}(y-x)\ln \left( \frac{y}{x}\right) \qquad \forall \kappa>1.
	\]
	
	Replacing $x$ by $\frac{ \displaystyle \prod\limits_{j=1}^{m-1} (a^n_j)^{\alpha_j}}{\phi^n}$ and $y$ by $\frac{a^n_m}{\phi^n}$, we obtain the following inequality
	
	\[
   \frac{ \displaystyle \prod\limits_{j=1}^{m-1} (a^n_j)^{\alpha_j}}{\phi^n} \leq \kappa \frac{a^n_m}{\phi^n}+\frac{1}{\ln{\kappa}}\frac{ \Big(a^n_m-\displaystyle \prod\limits_{j=1}^{m-1} (a^n_j)^{\alpha_j}\Big)}{\phi^n}\ln\left({\frac{a^n_m}{\prod\limits_{j=1}^{m-1} (a^n_j)^{\alpha_j}}}\right)
   \]
	
	The following two bounds one can deduce as a conclusion of the above algebraic inequality and relation \eqref{Uniform-integrability traingular}
	
	\begin{equation}\label{singular int}
		\left \{
		\begin{aligned}
			\bigint_{\Omega_T} \frac{ \displaystyle \prod\limits_{j=1}^{m-1} (a^n_j)^{\alpha_j}}{\phi^n} \leq 2TM_2+\frac{1}{\ln{2}}E(a_{i,0}:i=1,2\cdots,m) . \\
			\Bigg\{\frac{ \displaystyle \prod\limits_{j=1}^{m-1} (a^n_j)^{\alpha_j}}{\phi^n}\Bigg\}\bigg|_{n\in\mathbb{N}} \text{is uniformly integrable.}
		\end{aligned}
		\right .
	\end{equation}
	
	\vspace{.1cm}
	
	Define $g^n=\frac{\Big(a^n_m-\displaystyle \prod\limits_{j=1}^{m-1} (a^n_j)^{\alpha_j}\Big)}{\phi^n}$, the above bounds yields
	
	\begin{equation}\label{g_n L^1 norm}
		\left \{
		\begin{aligned}
			\Vert g^n \Vert_{L^1(\Omega_T)}\leq & M_2T+2TM_2+\frac{1}{\ln{2}}E(a_{i,0}:i=1,2\cdots,m). \\
			& \{g^n\}|_{n\in\mathbb{N}} \ \text{is \ uniformly \ integrable. }
		\end{aligned}
		\right .
	\end{equation}
	
	By Dunford Petties' theorem  $g^n$  is weakly compact. So from the property of heat kernel we can conclude  $\exists \ a_i$, such that,   $a_i^n$  converges to $a_i$   a.e in  $[0,T)\times {\Omega},  \forall i\in \Lambda_2$.  Vitali convergence lemma furthermore yields   $a_i^n \rightarrow a_i$  in  $L^1(\Omega_T),  \forall \ i\in \Lambda_2$.  
	
	\vspace{.1cm}
	
	Next we are moving to show uniform $L^{\infty}$ time space bound for all the species by showing some uniform $L^p(\Omega_T)$  norm bound, where $p$ can be large. All these bounds will be uniform with respect to the index $n$. The proposition is as follows

	\vspace{.3cm}

	\begin{Prop}\label{L^p estimation d_m, d_m>0}
		 There exists a positive constant $\Tilde{C}$,independent of the index '$n$' and depends only on initial condition and time $T$, such that
		\begin{align*}
			\Vert a^n_i \Vert_{L^{\infty}(\Omega_T)}< \Tilde{C}, \qquad \forall i=1,\cdots,m.
		\end{align*}
	\end{Prop}
	
	\vspace{.2cm}
	
	Proof: We have the following identity for \ $i \in \Lambda_1$
	
	\begin{gather}
		-\partial_t a^n_i =\partial_t a^n_m -d_m \Delta a^n_m, \qquad n.\nabla_{x}a^n_i=n.\nabla_{x}a^n_m=0.
		\label{equation relation}
	\end{gather}

	Consider   $\Theta \geq 0 \in C_c^{\infty}(\Omega_{\tau})$(space of all compactly supported smooth function),satisfies
	
	\begin{flalign*}
		-[\partial_t \phi+d_m \Delta\phi]=& \Theta \qquad \ \ \ \mbox{in}\ \Omega_{\tau}\\
		n.\nabla_{x} \phi=& 0 \qquad \mbox{on}\ (0,\tau)\times \partial\Omega\\
		\phi(\tau)=&0 \qquad \ \ \ \ \mbox{in}\ \Omega.
	\end{flalign*}
	
	We have  $\phi \geq 0$  and the following estimate \cite{Pie2010}\cite{quittner2019superlinear}\cite{Lam87}
	
	\[
	\Vert\phi_t \Vert_{L^q(\Omega_{\tau})}+\Vert\Delta \phi \Vert_{L^q(\Omega_{\tau})}+\sup\limits_{s\in[0,\tau]}\Vert\phi(s)\Vert_{L^q(\Omega)}+\Vert\phi\Vert_{L^q((0,\tau)\times \partial\Omega)}\leq C_{q,T,d_m}\Vert \Theta \Vert_{L^q(\Omega_{\tau})},
	\]
	
	where  $C_{q,T,d_m}>0$ is a constant and $q\in(1,\infty)$, arbitrary.
	
	\vspace{.2cm}

	Multiply equation \eqref{equation relation} by  $\phi$  and integration by parts yields
	
	\begin{flalign*}
		\int_{\Omega_{\tau}}a^n_m \Theta= \int_{\Omega}(a_{m,0}+a_{i,0})\phi(0)+\int_{\Omega_{\tau}}a^n_i \partial_t \phi.
	\end{flalign*}

	Let $p\in(1,\infty)$ be the H\'older conjugate of $q$. Apllying H\"older inequality on the previous relation we obtain
	
	\[ 
	\int_{\Omega_{\tau}}a^n_m \Theta\leq \big( \Vert a_{m,0}\Vert_{L^p(
		\Omega)}+\Vert a_{i,0}\Vert_{L^p(
		\Omega)}\big)\Vert \phi_{0}\Vert_{L^q(
		\Omega)}+\Vert a_i^n\Vert_{L^p(
		\Omega_{\tau})}\Vert \partial_{t} \phi\Vert_{L^q(
		\Omega_{\tau})}.
	\]

	Let  $C_{1}=C_{q,T,d_m}\max\{1,\Vert a_{m,0}\Vert_{L^p(
		\Omega)}+\Vert a_{i,0}\Vert_{L^p(\Omega)}\}$, then
	
	\[ 
	\int_{\Omega_{\tau}}a^n_m \Theta\leq C_{1}\big(1+\Vert a^n_i \Vert_{L^p(\Omega_{\tau})}\big) \Vert \Theta \Vert_{L^q(\Omega_{\tau})}.
	\]
	
	Hence by duality
	\begin{gather*} 
		\Vert a^n_m \Vert_{L^p(\Omega_{\tau})} \leq C_{1}\big(1+\Vert a^n_i \Vert_{L^p(\Omega_{\tau})}\big) \qquad \forall \ \tau\in(0,T].
	\end{gather*}
	
	Replacing $\tau$ by $t$, we can rewrite the relation as
	\begin{gather}
		\Vert a^n_m \Vert_{L^p(\Omega_{t})} \leq C_{1}\big(1+\Vert a^n_i \Vert_{L^p(\Omega_{t})}\big) \qquad \forall \ t
		\in(0,T]\label{relation 1 use}.
	\end{gather}

	Consider the o.d.e.
	\[
	  \partial_t a^n_i \leq a^n_m.
	\]
	
	 An application of Minkwoski's integral inequality and Jensen inequality yields,
	\[ 
	\Vert a^n_i \Vert^p_{L^p(\Omega)} \leq  2^{p-1}\Vert a_{i,0} \Vert^p_{L^p(\Omega_T)}+ 2^{p-1}T^{p-1}\int_{0}^{t}\Vert a^n_m \Vert^p_{L^p(\Omega)}.
	\]
	
	Raising  to the power $p$ in relation \eqref{relation 1 use}, yields
	
	\begin{gather*}
		\int_{0}^{t}\Vert a^n_m \Vert^p_{L^p(\Omega)}\leq 2^{(p-1)} C_{1}^p+2^{(p-1)} C_{1}^p\Vert a^n_i \Vert^p_{L^p(\Omega_t)}.
	\end{gather*}

	Consider the constant  $C_{2}=\max\{(1+T)^{p-1}2^{(2p-2)} C_{1}^p,2^{p-1}\Vert a_{i,0} \Vert^p_{L^p(\Omega_T)}:\forall i=1,\cdots,m \}$. From the two relations above we can conclude
	
	\[ 
	\Vert a^n_i \Vert^p_{L^p(\Omega)} \leq C_{2} \Big( 1+ \int_{0}^t \Vert a^n_i \Vert^p_{L^p(\Omega)} \Big).
	\]
	
	 Gr\"ownwall inequality yields

	\[ 
	\sup\limits_{t\in[0,T)}\Vert a^n_i \Vert_{L^p(\Omega)}\leq \left(1+e^{C_{2}T}+\max\limits_{\{i=1,\cdots,m\}}\Vert a_{i,0} \Vert_{L^p(\Omega)}e^{C_{2}T}\right)^{\frac{1}{p}},
	\]
	
	which implies
	
	\[
	\Vert a^n_i \Vert_{L^p(\Omega_T)}\leq \left(T+Te^{C_{2}T}+T\max\limits_{\{i=1,\cdots,m\}}\Vert a_{i,0} \Vert_{L^p(\Omega)}e^{C_{2}T}\right )^{\frac{1}{p}}.
	\]

	Relation \eqref{relation 1 use} further yields
	
	\[ 
	  \Vert a^n_m \Vert_{L^p(\Omega_T)}\leq C_1\left( 1+\left(T+Te^{C_{2}T}+T\max\limits_{\{i=1,\cdots,m\}}\Vert a_{i,0} \Vert_{L^p(\Omega)}e^{C_{2}T}\right )^{\frac{1}{p}}\right).
	\]

	Let's now consider the diffusing i.e., non-degenerate species.	For $i\in\Lambda_2$, we have the following relation
	
	\[ 
	\partial_t (a^n_i+a^n_m)-\Delta(d_i a^n_i +d_m a^n_m)=0, \qquad n.\nabla_{x}a^n_i=n.\nabla_{x}a^n_m=0.
	\]
	
	Here we use duality argument as in   theorem \ref{duality estimate}[see appendix].  There exists constant $C_{3}>0$, depends only on the initial condition, domain and the time $T$, such that the following relation holds
	
	\begin{gather}
		\Vert a^n_i \Vert_{L^p(\Omega_T)}\leq C_{3}, \qquad \forall \ i\in \Lambda_2, \ \forall n\in \mathbb{N}. \label{towards infinite estimation}
	\end{gather}
	
	Consider the following constant \begin{equation*}
		\hat{C}=\max \left \{
		\begin{aligned}
			&C_3,\\
			& \left ( 2^{(p-1)} C_{1}^p+2^{(p-1)} TC_{1}^p \left(C_2+C_2^2+C_{2}^2e^{C_{2}T}+C_2\max\limits_{\{i=1,\cdots,m\}}\Vert a_{i,0} \Vert_{L^p(\Omega)}e^{C_{2}T}\right)\right)^{\frac{1}{p}},\\
			& \left(TC_2+TC_2^2+TC_{2}^2e^{C_{2}T}+TC_2\max\limits_{\{i=1,\cdots,m\}}\Vert a_{i,0} \Vert_{L^p(\Omega)}e^{C_{2}T}\right )^{\frac{1}{p}}.
		\end{aligned}
		\right .
	\end{equation*}

The following norm bound holds
	
	\[
	\Vert a^n_i \Vert_{L^p(\Omega_T)}\leq \hat{C}, \quad \forall i=1,\cdots,m.
	\]
	
	If we in particular take  $p=2NQ$  then we have $g^n=\frac{\Big(a^n_m-\displaystyle \prod\limits_{j=1}^{m-1} (a^n_j)^{\alpha_j}\Big)}{\phi^n} \in L^{2N}(\Omega_T)$. The $L^{\infty}$ time space bound is a consequence of theorem \ref{estimation 1}[see appendix]. For a positive constant $\Tilde{C}>0$, independent of the index $n$, we have the following bound

	\begin{gather}
		\Vert a^n_i \Vert_{L^{\infty}([0,T)\times\Omega)}\leq \Tilde{C}, \qquad \forall  i\in 1,\cdots,m, \  \forall \ n\in \mathbb{N}. \label{ infinite estimation} 
	\end{gather}
	
	\vspace{.2cm}
	
	In the next proposition we find a weak limit for the degenerate species.

	\begin{Prop}
	For all $i\in \Lambda_1$, $a_i^n$ has an almost everywhere convergent subsequence, converges to $a_i$, which weakly satisfies the equation
		\[
		\partial_t a_i =a_m-\displaystyle\prod\limits_{j=1}^{m-1}(a_j)^{\alpha_j}.
		\]
	\end{Prop}
	
	\vspace{.1cm}
	
	Proof: Let  $i,j\in\Lambda_1$, then  $a^n_i-a^n_j=a_{i,0}-a_{j,0}=\phi_{i,j}$.  We fix a  $j\in\Lambda_1$, such that $\alpha_j\in \{1\}\cup[2,\infty)$\eqref{SC}.
	
	\[
	\partial_t a^n_j = \frac{a^n_m-\displaystyle\prod\limits_{j=1}^{m-1}(a^n_i)^{\alpha_i}}{\phi^n}.
	\]
	
	Let us define some quantities
	
	\begin{flalign*}
		\delta^1_n=&(\phi_n)^{-1}, \\
		\delta^2_n(a^n_j)=& {(a^n_j)}^{\alpha_j-1} \displaystyle\prod\limits_{i\in \Lambda_1}(\phi_{ij}+a^n_j)^{\alpha_i},\\
		\delta_n^3=&\displaystyle\prod\limits_{i\in \Lambda_2\setminus \{m\}}{(a^n_i)}^{\alpha_i}.\\
	\end{flalign*}
	
	We can rewrite the above relation as follows
	
	\begin{gather} 
		a^n_j(t,x)=a_{j,0}(x)e^{-\int_{0}^{t}\delta^n_1 \delta^n_2(a^n_j) \delta^n_3ds}+\int_{0}^{t}a_m^n \delta^n_1e^{-\int_{s}^{t}\delta^n_1 \delta^n_2(a^n_j) \delta^n_3 d\sigma}ds, \label{relation 1}
	\end{gather} 
	
	Thanks to proposition \eqref{L^p estimation d_m, d_m>0} we obtain the following pointwise bound
	
	\[
	   a^n_j(t,x) \leq a_{j,0}(x)+T\Tilde{C} \quad \text{for} \ x\in\Omega \ \text{a.e.}
	\]
	
	We can conclude  $a^n_j(t,x)$  bounded uniformly for a given  $x\in \Omega$ a.e., which implies  $\delta_n^1\rightarrow 1$\ pointwise a.e. Furthermore Vitali convergence lemma yields $\delta_n^1\rightarrow 1 \ \text{in}\ L^1(0,T)$ for a given  $x\in \Omega$  a.e.
	
	Let's introduce a Picard type iteration on \eqref{relation 1}
	
	\begin{flalign*}
		a^n_{j,p+1}(t,x)=a_{j,0}(x)&e^{-\int_{0}^{t}\delta^n_1 \delta^n_2(a^n_{j,p}) \delta^n_3ds}+\int_{0}^{t}a_m^n \delta^n_1e^{-\int_{s}^{t}\delta^n_1 \delta^n_2(a^n_{j,p}) \delta^n_3 d\sigma}ds,\\
		a^n_{j,p}(0,x)=& a_{j,0} \qquad \forall  n\in\mathbb{N}.
	\end{flalign*}

	Again thanks to proposition \eqref{L^p estimation d_m, d_m>0}, we obtain the following pointwise bound
	
	\[ 
	\sup\limits_{t\in[0,T]}\sup\limits_{p\in\mathbb{N}}\vert a^n_{j,p}(t,x)\vert \leq a_{j,0}(x)+T\Tilde{C}=C_{4} <+\infty, \quad \text{for}\ x\in\Omega \ \text{a.e.}
	\]

	\begin{flalign*}
		\vert  & a^n_{j,p+1}(t,x)- a^n_{j}(t,x) \vert \leq\\
		& \Bigg((1+T)\Big(a_{j,0}(x)+T\Tilde{C}\Big)\displaystyle\prod\limits_{i\in \Lambda_2\setminus \{m\}}\Tilde{C}^{\alpha_i} \sup\limits_{r\in[0,C_{4}]}\Big [\Big \vert \frac{d(\zeta(r))}{dr}\Big \vert \Big]\Bigg )\int_{0}^{t}\vert  a^n_{j,p}(s,x)- a^n_{j}(s,x) \vert ds,
	\end{flalign*}

	\[ 
	\text{where}\ \zeta(r)=r^{\alpha_j-1}\displaystyle\prod\limits_{i\in\Lambda_1}(\phi_{i,j}+r)^{\alpha_i}.
	\]

	Consider the constant $C_{5}=(1+T)\Big(a_{j,0}(x)+T\Tilde{C}\Big)\displaystyle\prod\limits_{i\in \Lambda_2\setminus \{m\}}\Tilde{C}^{\alpha_i} \sup\limits_{r\in[0,C_{4}]}\Big [\Big \vert \frac{d(\zeta(r))}{dr}\Big \vert \Big]$. The previous relation can be rewritten as

	\[ 
	\vert  a^n_{j,p+1}(t,x)- a^n_{j}(t,x) \vert \leq C_{5}\int_{0}^{t}\vert  a^n_{j,p}(s,x)- a^n_{j}(s,x) \vert ds.
	\]

	Induction on $p$ yields  
	
	\[ 
	\vert  a^n_{j,p+1}(t,x)- a^n_{j}(t,x) \vert \leq \frac{C_{5}^{p} t^{p}}{
		p!}\int_{0}^{t}\vert a^n_{j,0}(s,x)- a^n_{j}(s,x) \vert \leq 2TC_{4}\frac{C_{5}^{p} t^{p}}{p!}, \quad \forall  x\in\Omega  \ \text{a.e.}
	\]
	
	\[ 
	\sup\limits_{n\in\mathbb{N}}\sup\limits_{t\in[0,T]}\vert  a^n_{j,p+1}(t,x)- a^n_{j}(t,x) \vert \leq 2TC_{4}\frac{C_{5}^{p} T^{p}}{p!}.
	\]
	
	From the above relation we can conclude
	
	\begin{gather}
		\lim\limits_{p\rightarrow\infty}\sup\limits_{n\in\mathbb{N}}\sup\limits_{t\in[0,T]}\vert  a^n_{j,p}(t,x)- a^n_{j}(t,x) \vert=0 \label{relation 2}.
	\end{gather}

	For $p\geq1$, a similar calculation on  $a^n_{j,p}(t,x)$  yields the following estimates

	\[ 
	\vert  a^n_{j,p+2}(t,x)- a^n_{j,p+1}(t,x) \vert \leq \frac{C_{5}^{p-1} t^{p-1}}{
		(p-1)!}\int_{0}^{t}\vert a^n_{j,1}(s,x)- a^n_{j,0}(s,x) \vert \leq 2TC_{4}\frac{C_{5}^{p-1} t^{p-1}}{(p-1)!}.
	\]

	\[ 
	 \implies\sup\limits_{n\in\mathbb{N}}\sup\limits_{t\in[0,T]}\vert  a^n_{j,p+2}(t,x)- a^n_{j,p+1}(t,x) \vert \leq 2TC_{4}\frac{C_{5}^{p-1} T^{p-1}}{(p-1)!}.
	\]
	
		From the above relation we obtain the following limit
	
	\begin{gather}
		 \lim\limits_{p\rightarrow\infty}\sup\limits_{n\in\mathbb{N}}\sup\limits_{t\in[0,T]}\vert  a^n_{j,p+r}(t,x)- a^n_{j,p}(t,x).\vert=0, \qquad \forall r\in \mathbb{N}.\label{relation 3}.
	\end{gather}

	\vspace{.3cm}
	
	 Vitali convergence lemma yields the two following limits 
	
	\begin{flalign}
		\lim\limits_{n\rightarrow\infty}\int_{0}^{t}\vert \displaystyle\prod\limits_{i\in \Lambda_2\setminus \{m\}}a_i^{\alpha_i}-\delta_1^n \delta_3^n \vert = 0 \qquad \forall x\in \Omega \ \text{a.e.,}\ \forall  t \in [0,T), \label{relation 4}\\
		\lim\limits_{n\rightarrow\infty}\int_{0}^{t}\vert a_m-a_m^n\delta^n_1\vert=0 \qquad \forall x\in \Omega \ \text{a.e.,}\ \forall  t \in [0,T). \label{relation 5}
	\end{flalign}

	for a given  $x\in \Omega$ a.e., let us consider another Picard type iteration

	\begin{flalign*}
		b_{p+1}(t,x)=a_{j,0}e^{-\int_0^t \delta^2_n(b_p)\delta^4_n}&+\int_{0}^t a_me^{-\int_s^t \delta^2_n(b_p)\delta^4_n d\sigma}ds,\\
		b_{0}=a_{j,0},
	\end{flalign*}

	\[
	\text{where}\ \ \delta_n^4=\displaystyle\prod\limits_{i\in\Lambda_2\setminus \{m\}}a_i^{\alpha_i}.
	\]
	
	Thanks to proposition \eqref{L^p estimation d_m, d_m>0} we obtain the following pointwise bound

	\[ 
	\sup\limits_{t\in[0,T]}\sup\limits_{p\in\mathbb{N}}\vert b_p(t,x)\vert \leq a_{j,0}(x)+T\Tilde{C}=C_{4} <+\infty, \quad x\in\Omega \ \text{a.e.}
	\]

	For  $p=0$,  we have  $\lim\limits_{n\rightarrow\infty}\sup\limits_{t\in[0,T]}\vert a^n_{j,0}-b_0\vert=0$. Induction on $p$ yields

	\begin{flalign*}
		\vert a^n_{j,p+1}-b_{p+1}\vert \leq & a_{j,0}\Big[ \Big \vert \int_{0}^{t}(\delta_n^4-\delta^1_n\delta^3_n)\delta^2_n(a^n_{j,p}) \Big\vert+\sup\limits_{r\in[0,\Tilde{C}+C_{4}]}\Big [\Big \vert \frac{d(\zeta(r))}{dr}\Big \vert \Big]\int_0^t \vert a^n_{j,p}-b_p \vert \Big]
		\\
		+ \int_{0}^{t}\vert a_m & -a_m^n\delta^n_1 \vert e^{-\int_s^t \delta^1_n \delta^2_n \delta^3_n}+(C_{4}+\Tilde{C})T\sup\limits_{r\in[0,\Tilde{C}+C_{4}]}\Big [\Big \vert \frac{d(\zeta(r))}{dr}\Big \vert \Big]\int_0^t \vert a^n_{j,p}-b_p \vert 
		\\
		+& (\Tilde{C}+C_{4})T\Big\vert \int_0^t (\delta^4_n-\delta^1_n\delta^3_n)\delta^2_n(a^n_{j,p}) \Big\vert.
	\end{flalign*}

	Consider $\lim\limits_{n\rightarrow\infty}\sup\limits_{t\in[0,T]}\vert a^n_{j,p}-b_p\vert=0$. Then Vitali convergence lemma along with the relations \eqref{relation 4} and \eqref{relation 5} implies 
	
	\[
	\lim\limits_{n\rightarrow\infty}\sup\limits_{t\in[0,T]}\vert a^n_{j,p+1}-b_{p+1}\vert=0.
	\] 
	
	So by induction on  $p$, we obtain

	\begin{gather}
		\lim\limits_{n\rightarrow\infty}\sup\limits_{t\in[0,T]}\vert a^n_{j,p}-b_p\vert=0 \qquad \forall p\in \mathbb{N}, \  \forall x\in \Omega \ \mbox{a.e.}\label{relation 6}
	\end{gather}

	Triangle inequality yields
	
	\[ 
	\vert b_{p+r}-b_p\vert \leq \vert a^n_{j,p+r}-b_{p+r}\vert+\vert a^n_{j,p}-b_p\vert+\vert a^n_{j,p+r}-a^n_{j,p}\vert, \quad  \forall \ r\in\mathbb{N}.
	\]

	From the two relations\eqref{relation 6}  and \eqref{relation 3}, we conclude

	\[
	\lim\limits_{n\rightarrow\infty}\sup\limits_{t\in[0,T]}\vert b_{p+r}-b_p\vert=0 \qquad \forall \ r\in\mathbb{N},
	\]
	
	which further implies
	
	\[
	 \{b_p(t)\} \ \text{is \ Cauchy\ in}\ L^{\infty}[0,T]\ \ \text{for\ } x\in\Omega \ \text{a.e.}
	\]

	Let  $b_p(t,x)$ converges to  $a_j(t,x)  $  for a given  $x\in\Omega$  a.e, $\forall t\in[0,T]$. Furthermore as this  convergence is uniform with respect to time we can take the limit inside the integral. Hence for all $x\in\Omega$ a.e., $a_j(t,x)$ satisfies the following relation

	\[
	a_j(t,x)=a_{j,0}(x)e^{-\bigint_{0}^{t}\displaystyle\prod_{1}^{m-1} (a_j)^{\alpha_j}}+\large{\int^{t}_{0}} a_me^{-\bigint_{s}^{t}\displaystyle\prod_{1}^{m-1} (a_j)^{\alpha_j}d\sigma}ds \qquad \forall t\in[0,T]
	\]

	\vspace{.2cm}
	
	From the relations \eqref{relation 2} and  \eqref{relation 3}, for a given  $x\in\Omega$  a.e.,  we can conclude   $a^n_{j}(t,x)$ converges to  $a_j(t,x)$  uniformly  $\forall t\in[0,T]$.   Now this pointwise limit and uniform integrability character of  $a^n_{j}(t,x)$,(from relation \eqref{Uniform-integrability traingular}) provides us (Vitali converges lemma)

	\begin{gather}
		a^n_{j}(t,x) \rightarrow a_j(t,x) \ \text{in} \ \
		L^1(\Omega_T) \label{picard method}.
	\end{gather}

	As  $a^n_i-a^n_j=a_{i,0}-a_{j,0}=\phi_{i,j},\  \forall \ i,j\in\Lambda_1$,  we have

	\[ 
	a^n_{i}(t,x) \rightarrow a_i(t,x)=\phi_{i,j}-a_j(t,x) \ \text{in} \ \ L^1(\Omega_T), \ \ \ \forall i\in\Lambda_1. 
	\]

	\vspace{.3cm}

		Proof of theorem \ref{Existence A1}:
		The weak limit  $(a_1,a_2,\cdots,a_m)$  is a solution of our triangular system-$\text{A}_1$(in distributional sense), where  $a_i$  satisfies the following pointwise bound(proposition \eqref{L^p estimation d_m, d_m>0})
		
		\[
		\Vert a_i \Vert_{L^{\infty}(\Omega_T)} \leq \Tilde{C} \qquad \forall \ i=1,2,\cdots,m.
		\]
		
		Therefore parabolic regularity gives us a smooth unique solution
		\cite{DF15} \cite{EMT20} \cite{FMT20} \cite{quittner2019superlinear}.

	\section {\textbf{\large{Existence of solution of degenerate triangular system-$\textbf{A}_2$}}} 
	
	We need uniform $L^{\infty}$ time-space  bound of the solution of approximate solution, where the bound will be uniform with respect to index $n$. It can be derived easily from the degenerate nature of system-$\text{A}_2$, as described in the following proof.
	
	\vspace{.2cm}
	
	Proof of theorem-\ref{Existence A2}: Just as in the previous section we consider the approximate system \eqref{approximate  system} with $d_m=0$ and one or more  $d_i=0$, where  $i\in\{1,2,\cdots,m-1\}$. Addition of a degenerate equation corresponding to $a_i$ where $i\neq m$ with that of $a_m$, yields
		
		\[ 
		a^n_m+a^n_i=a_{m,0}+a_{i,0} \ \ \text{i.e.,} \ \  \ a^n_m \ \  \text{uniformly \ bounded}.
		\]
		
		Similar duality arguments as in the previous section yields
		
		\[
		 a_j^n\ \text{all\ uniformly\ bounded},\ \ \forall j=1,2,\cdots,m.(\text{theorem}\ref{estimation 1}[\text{see appendix}]).
		\]
		
		Hence we have existence of unique smooth global   solution.

	\section {\textbf{\large{Existence of solution of degenerate triangular system-$\textbf{A}_3$}}} 
	
	It is evident from our discussion in first section that the rate function for the approximate system $g^n=\frac{\Big(a^n_m-\displaystyle \prod_{1}^{m-1} (a^n_j)^{\alpha_j}\Big)}{\phi^n}$ are uniformly bounded in $L^1(\Omega_T)$ and uniformly integrable\eqref{g_n L^1 norm}. Furthermore the approximate solutions $a^n_i$ are also uniformly bounded in $L^1(\Omega_T)$ and uniformly integrable\eqref{Uniform-integrability traingular}. Based on these we can extract a weak solution for the system-$\text{A}_3$.
	
	\vspace{.3cm}
	
Proof of theorem-\ref{Existence A3}, First part, weak solution: \ By Dunford Petties' theorem  $g^n$ is weakly compact and   $\exists a_i$  such that   $a_i^n$  converges to $a_i$   a.e in  $[0,T)\times\Omega$, thanks to the property of heat  kernel, $\forall i=1,2,\cdots,m-1$. Hence by Vitali convergence lemma  $a_i^n \rightarrow a_i$  in  $L^1(\Omega_T),  \forall i=1,2,\cdots,m-1$.
		The differntial equation corresponds to $a_m^n$ is the following
		
		\[ 
		\partial_t a^n_m= \frac{\displaystyle\prod\limits_{j=1}^{m-1}(a^n_j)^{\alpha_j}-a^n_m}{\phi^n}.
		\]

		\[
		a^n_m=\bigint_{0}^{t}\frac{\displaystyle\prod\limits_{j=1}^{m-1}(a^n_j)^{\alpha_j}}{\phi^n} e^{(-\int_{s}^{t} \frac{1}{\phi^n}d\sigma)} ds.
		\]
		
		From relation \eqref{singular int}, we have 
		
		\[
		\bigint_{\Omega_T}\frac{ \displaystyle \prod\limits_{j=1}^{m-1} (a^n_i)^{\alpha_j}}{\phi^n} \leq 2TM_2+\frac{1}{\ln{2}}E(a_{i,0}:i=1,2\cdots,m),
		\]
		 which implies  $a_m^n(t,x)$  bounded for  $\forall t\in(0,T),x\in {\Omega}$ a.e. So  $\phi^n$  converges to 1 pointwise a.e. \\
		Uniform integrability of   $\Bigg\{\frac{ \displaystyle \prod\limits_{j=1}^{m-1} (a^n_j)^{\alpha_j}}{\phi^n}\Bigg\}\Bigg|_{n\in\mathbb{N}}$(relation-\eqref{singular int}) implies  $\frac{ \displaystyle \prod\limits_{j=1}^{m-1} (a^n_j)^{\alpha_j}}{\phi^n}\xrightarrow[n\rightarrow\infty] \displaystyle\prod\limits_{j=1}^{m-1}a_j^{\alpha_j}$ in  $L^1(\Omega_T)$(by Vitali convergence lemma). Furthermore Fatou's lemma yields

		\[ 
		0\leq\int_{\Omega}\liminf{\int_{0}^T} \Bigg\vert\frac{ \displaystyle\prod\limits_{j=1}^{m-1}(a^n_j)^{\alpha_j}}{\phi^n} -\displaystyle\prod\limits_{j=1}^{m-1}(a_j)^{\alpha_j} \Bigg\vert \leq \liminf \int_{\Omega}\int_{0}^T \Bigg\vert\frac{ \displaystyle\prod\limits_{j=1}^{m-1}(a^n_j)^{\alpha_j}}{\phi^n} -\displaystyle\prod\limits_{j=1}^{m-1}(a_j)^{\alpha_j} \Bigg\vert\leq0.
		\]
		
		It yields the following limit
		
		\[
		 \liminf{\int_{0}^T} \Bigg\vert\frac{ \displaystyle\prod\limits_{j=1}^{m-1}(a^n_j)^{\alpha_j}}{\phi^n} -\displaystyle\prod\limits_{j=1}^{m-1}(a_j)^{\alpha_j} \Bigg\vert=0.
		\]

		We can conclude from the above relation that there exists a subsequence indexed by  $n_k$ ( without loss of generality we take $n_k$ as  $n$), such that for all  $x\in {\Omega}$ a.e.,
		
		\[
		\lim \bigint_{0}^{t}\frac{\displaystyle\prod\limits_{j=1}^{m-1}(a^n_j)^{\alpha_j}}{\phi^n} e^{(-\int_{s}^{t} \frac{1}{\phi^n}d\sigma)} ds=\bigint_{0}^{t}\displaystyle\prod\limits_{j=1}^{m-1}(a_j)^{\alpha_j}e^{(s-t)}ds=a_m.
		\]
		
		Furthermore  uniform integrability character  \eqref{Uniform-integrability traingular} yields  $a^n_m$ converges to  $a_m$  in  $L^1(\Omega_T)$( Vitali converges lemma).
		
		\vspace{.2cm}

		Consider the following relation
		
		\[
		\partial_t (a^n_i+a^n_m)-d_i\Delta a^n_i=0, \qquad n.\nabla_{x}a^n_i=n.\nabla_{x}a^n_m=0
		\]
		
		The following uniform  integral estimate holds (see appendix theorem\ref{2nd order estimate})

		\begin{flalign} 
			\int_0^{T} \int_{\Omega}(a_i^n)^2+ a^n_ia^n_m \leq C_{\Omega}(1+T),  \qquad \forall T\geq 0, \quad \forall i=1,2,\cdots,m-1, \label{towards gradient estimate}
		\end{flalign}
	where the constant $C_{\Omega}>0$, independent of index $n$.

		 Next we multiply $a_i^n$ to the equation corresponding to $a_i^n$ in the approximate system and integrate over  $\Omega_T$. We obtain the following relation

		\begin{align}\label{L^2 gradient estimate}
			\frac{1}{2}\int_{\Omega}(a_i^n)^2+\int_{\Omega_T}\vert \nabla a_i^n\vert^2 \leq \frac{1}{2}\int_{\Omega}(a_{i,0})^2+\int_{\Omega_T}\frac{a^n_ia^n_m}{\phi^n}, \qquad \forall  i=1,2,\cdots,m-1.
		\end{align}

		Above relation yields a $n$ index free uniform bound for $ \nabla a_i^n$ in  $L^2(\Omega_T)$, for all $i=1,\cdots,m-1$. So it  weakly converges (upto a subsequence) to   $\nabla a_i$(in distributional sense), for all  $i=1,2,\cdots,m-1$. We can conclude $(a_1,a_2,\cdots,a_m)$ is a global in time positive weak solution of triangular system-$\text{A}_3$\eqref{A3}($d_m=0$).

	\vspace{.3cm}
	
	We can extend this weak solution to classical solution for dimension $N=1,2$,  by the $L^p$ integral estimation of Neumann green function for the heat equation.

	\vspace{.3cm}

	Proof of theorem-\ref{Existence A3}, second part, smooth solution:  $a_i$ satisfies the following equation weakly
		
		\begin{equation} \nonumber
			\left\{
			\begin{aligned}
				\partial_t a_i-d_i \Delta a_i = & \psi=a_m-\displaystyle{\prod\limits_{j=1}^{m-1}a_j^{\alpha_j}} \qquad \qquad \qquad \ \ \  \mbox{ in } \Omega_T,\forall i=1,\cdots,m-1\\
				\nabla a_i .\gamma=& 0 \qquad \qquad \qquad \qquad \qquad \ \ \ \mbox { on } (0,T)\times\partial\Omega\\
				a_i(0,x)=& a_{i,0} \qquad \qquad \qquad \qquad\qquad \qquad \quad \ \mbox{ in }\Omega.
			\end{aligned}
			\right.
		\end{equation}
		
		Let  $G_i(t,s,x,y)$ be the Neumann Green function corresponding to the operator $\partial_t-d_i \Delta$, we can  express $a_i$ in the following way
		
		\begin{align}\label{d_m: representation of solution}
			a_i(t,x)= \Tilde{a}_i(t,x)+\int_{0}^{t}G_{i}(t,s,x,y)\psi(s,y)dyds,
		\end{align}
		
		where $\Tilde{a}_i(t,x)$ satisfies the following pde,
		
		\begin{equation} \nonumber
			\left\{
			\begin{aligned}
				\partial_t \Tilde{a}_i-d_i \Delta a_i = & 0 \qquad \qquad \qquad \ \ \  \mbox{ in } \Omega_T,\forall i=1,\cdots,m-1\\
				\nabla \Tilde{a}_i .\gamma=& 0 \qquad \qquad \quad \ \ \ \mbox { on } (0,T)\times\partial\Omega\\
				\Tilde{a}_i(0,x)=& a_{i,0} \qquad \qquad \qquad \  \mbox{ in }\Omega.
			\end{aligned}
			\right.
		\end{equation}

		We use the following estimates from the articles \cite{Morra83}\cite{ML15}. There exists $\texttt{C}_{H},\kappa,C_{S}>0$,such that
		
		\begin{equation}\label{heat kernel estimate}
			\left\{
			\begin{aligned}
				0\leq  G_{i}(t,s,x,y)  \leq & \texttt{C}_{H} \frac{1}{(t-s)^{\frac{N}{2}}}e^{-\kappa\frac{\Vert x-y \Vert^2}{(t-s)}}=g(t-s,x-y), \quad \forall i=1,\cdots,m-1,\\
				\Vert  \Tilde{a}_i(t,x) \Vert_{L^p(\Omega)}\leq & C_{S}\Vert  a_{i,0} \Vert_{L^p(\Omega)},\qquad \forall i=1,\cdots,m-1.
			\end{aligned}
			\right.
		\end{equation}
		
		Non-negativity of solution yields
		
		\[
		a_i(t,x) \leq \Tilde{a}_i(t,x)+\int_{0}^{t}\int_{\Omega}G_{i}(t,s,x,y)a_m(s,y)dyds.
		\]
		
		Applying Minkwoski's integral inequality we obtain
		
		\[
		\Vert a_i(t,x) \Vert_{L^p(\Omega)} \leq C_{S,0}\Vert  a_{i,0} \Vert_{L^p(\Omega)}+ \int_{0}^{t} \Vert g_{j}(t-s,x)\Vert_{L^p(\Omega)} \Vert a_m \Vert_{L^1(\Omega)}.
		\]
		
		There exists a positive constant $C_{H,p}$, such that
		
		\begin{align} \label{d_m; Lp estimation a_i}
			\Vert a_i(t,x) \Vert_{L^p(\Omega)} \leq C_{S}\Vert  a_{i,0} \Vert_{L^p(\Omega)}+C_{H,p}M_2 \int_{0}^{t} (t-s)^{-\frac{N}{2}(1-\frac{1}{p})}.
		\end{align}
		
		In dimension $N=1,2$, the integral $\int_{0}^{t} (t-s)^{-\frac{N}{2}(1-\frac{1}{p})} < +\infty, \ \forall p\in[1,\infty)$. Hence the $p=4Q$ yields us the following relation
		
		\[
		\Vert a_i(t,x) \Vert_{L^{4Q}(\Omega)} \leq K_p, \qquad \forall i=1,\cdots,m-1, 
		\]
		where $K_p$ is a positive constant.
		
		Next consider the o.d.e.,
		 \[
		 \partial_{t}a_m \leq \displaystyle\prod\limits_{j=1}^{m-1}a_j^{\alpha_j}
		 \]
		 
		  Applying H\"older inequality on the integral relation corresponding to above differential equation, we obtain
		
		\[
		\Vert a_m(t,x) \Vert_{L^4(\Omega)} \leq \int_{0}^{t} \sum\limits_{i=1}^{m-1}\Vert a_i(t,x) \Vert_{L^{4Q}(\Omega)}^{Q} \leq mTK_{p}^{Q}.
		\]
		
		Above bounds yields $\psi \in L^{4}(\Omega_T)$. As $4>\frac{N+2}{2}$ for $N=1,2$, by  integrability estimation (see appendix theorem\ref{estimation 1})
		
		\[
		\Vert a_i(t,x) \Vert_{L^{\infty}(\Omega_T)} < +\infty, \qquad \forall i=1,\cdots,m.
		\]
		
		Next parabolic regularity gives us unique global smooth solution for dimension $N=1$ and $2$\cite{quittner2019superlinear}\cite{EMT20}.

	\vspace{.3cm}
	
	Next we discuss two special cases, one  where the rate function grows at most quadratically and another a 3 dimensional case. The three species degenerate triangular model as discussed in the article \cite{DF15} has  rate function with quadratic growth. We have the following proposition.
	
	\vspace{.3cm}

	\begin{Prop}\label{quadratic growth d_m=0}
		If $\displaystyle{G=\sum\limits_{i=1}^{m-1}\alpha_i\leq2}$, unique global in time classical solution exists upto dimension $N=5$ for system-$\text{A}_3$\eqref{A3}.
	\end{Prop}

	\vspace{.2cm}

	It is enough to show global in time classical solution exists for $G=2$ in dimension $N=3,4,5$. From representation of solution \eqref{d_m: representation of solution},  we have
	
	\[
	a_i(t,x) \leq \Tilde{a}_i(t,x)+\int_{0}^{t}\int_{\Omega}G_{i}(t-s,x,y)a_m(s,y)dyds.
	\]
	
	Let $1+\frac{1}{p}=\frac{1}{r}+\frac{1}{q}$. Applying Minkwoski's integral inequality and Young convolution inequality, we obtain
	
	\[
	\Vert a_i(t,x) \Vert_{L^p(\Omega)} \leq C_{S,0}\Vert  a_{i,0} \Vert_{L^p(\Omega)}+ \int_{0}^{t} \Vert g_{j}(t-s,x)\Vert_{L^{r}(\Omega)} \Vert a_m \Vert_{L^q(\Omega)}.
	\]
	
	There exists a positive constant $C_{H,r}$, such that\cite{ML15}\cite{Morra83}
	
	\begin{align}\label{Lp-Lq estimate}
		\Vert a_i(t,x) \Vert_{L^p(\Omega)} \leq C_{S,0}\Vert  a_{i,0} \Vert_{L^p(\Omega)}+ C_{H,r} \int_{0}^{t} t^{-\frac{N}{2}\big(\frac{1}{q}-\frac{1}{p}\big)}\Vert a_m \Vert_{L^q(\Omega)}.
	\end{align}
	
	 From Sobolev embedding theorem, we have  there exists positive constant $C_{SE}$, depending only on the domain,such that 
	 \[
	 \Vert a_i \Vert_{L^{\frac{2N}{N-2}}(\Omega)}\leq C_{SE}(\Vert \nabla a_i \Vert_{L^2(\Omega)}+ \Vert a_i \Vert_{L^2(\Omega)}), \qquad  \forall i=1\cdots,m-1.
	 \]
	 
	 Consider the o.d.e.,
	 \[
	 \partial_t a_m \leq \displaystyle\prod_{j=1}^{m-1} a_j^{\alpha_j}
	 \] 
	Applying Minkwoski's inequality and H\"older inequality on the integral relation corresponding to above differential equation, we obtain
	
	\begin{flalign*}
		\Vert a_m \Vert_{L^{\frac{N}{N-2}}(\Omega)} \leq \Vert a_{m,0} & \Vert_{L^{\frac{N}{N-2}}(\Omega)} + \int_{0}^{t} \sum\limits_{i=0}^{m-1} \Vert a_i \Vert_{L^{\frac{2N}{N-2}}(\Omega)}^2 \\
		& \leq \Vert a_{m,0} \Vert_{L^{\frac{N}{N-2}}(\Omega)} + \int_{0}^{t} \sum\limits_{i=0}^{m-1} 2C_{SE}^2(\Vert \nabla a_i \Vert_{L^{2}(\Omega)}^2+\Vert a_i \Vert_{L^2(\Omega)}^2).
	\end{flalign*}
	
	Relation \eqref{L^2 gradient estimate} and \eqref{towards gradient estimate} yields
	
	\[
	\Vert a_m \Vert_{L^{\frac{N}{N-2}}(\Omega)} \leq \Vert a_{m,0} \Vert_{L^{\frac{N}{N-2}}(\Omega)}+C_{SE}^2\sum\limits_{i=1}^{m-1}\Vert a_{i,0} \Vert_{L^2(\Omega)}+4mC_{SE}^2C_{\Omega}(1+T)=C_1.
	\]
	
	For dimension $N=3$, $a_m \in L^3(\Omega)$. In the relation \eqref{Lp-Lq estimate} we take $q=3$ and $p=+\infty$. It yields  
	
	\[
	\Vert a_i \Vert_{L^{\infty}(\Omega)} \leq C_{S}\Vert a_{i,0}\Vert_{L^{\infty}(\Omega)}+2C_{H,\frac{3}{2}}\sqrt{T}, \qquad \forall i=1,\cdots m-1.
	\]
	
This further implies $a_m \in L^{\infty}(\Omega_T)$. Hence smooth unique global in time solution exists in dimension $N=3$.
	
	Similarly for dimension $N=4$ and $5$, our goal will be to prove $a_i \in L^{\infty}(\Omega_T) \ \forall i=1,\cdots,m$. For this, it is enough to show that $a_m \in L^{q}(\Omega)$, where $q>\frac{N}{2}$. Rest follows from \eqref{Lp-Lq estimate}.
	
	\vspace{.3cm}
	
	\textbf{N=4 CASE:} We obtain $a_m \in L^{2}(\Omega)$. In relation \eqref{Lp-Lq estimate} we take $q=2$ and $p=6$. It yields $a_i \in L^{6}(\Omega)$. Let $a_i \in L^{2p}(\Omega)$, then application of Minkwoski's integral inequality and H\"older inequality provide us
	
	\begin{align}\label{L2p-Lp estimate}
		\Vert a_m \Vert_{L^p(\Omega)} \leq \Vert a_{m,0} \Vert_{L^p(\Omega)}+ \sum\limits_{i=1}^{m-1} \Vert a_i \Vert_{L^{2p}(\Omega_T)}^2,
	\end{align}
	
	which implies $ a_m  \in L^{p}(\Omega)$. As we have $a_i \in L^6(\Omega)$ for all $i =1,\cdots,m-1$, we obtain $a_m\in L^3(\Omega)$. As $3> \frac{4}{2}$, relation \eqref{Lp-Lq estimate} yields $a_i \in L^{\infty}(\Omega_T)$ for all $i=1,\cdots,m$. Hence smooth unique global in time solution exists.
	
	\vspace{.3cm}
	
	\textbf{N=5 CASE:} We obtain $a_m \in L^{\frac{5}{3}}(\Omega)$. Relation \eqref{Lp-Lq estimate} yields $a_i \in L^{\frac{5}{1.1}}(\Omega)$(taking $q=\frac{5}{3},p=\frac{5}{1.1})$. Relation \eqref{L2p-Lp estimate} further provides $a_m \in L^{\frac{5}{2.2}}(\Omega)$. Again we use use  relation \eqref{Lp-Lq estimate} with $q=\frac{5}{2.2}$ and $p=\frac{50}{3}$. We obtain $a_i \in L^{\frac{50}{3}}(\Omega)$. Relation \eqref{L2p-Lp estimate} further yields $a_m \in L^{\frac{25}{3}}(\Omega)$. As $\frac{25}{3}>\frac{5}{2}$, relation \eqref{Lp-Lq estimate} yields $a_i\in L^{\infty}(\Omega_T)$. Hence smooth unique global in time solution exists.

	\vspace{.3cm}
	
	Next we discuss a special 3 dimensional case. The proposition is as follows.
	
	\vspace{.2cm}
	
	\begin{Prop}
		Let dimension $N=3$. Then system-$\text{A}_3$\eqref{A3} has
		unique global in time classical solution for $\displaystyle{G=\sum\limits_{i=1}^{m-1}\alpha_i\leq \frac{10}{3}}$.
	\end{Prop}

	\vspace{.2cm}

	It is enough to show for $G=\frac{10}{3}$. In relation \eqref{Lp-Lq estimate} we take $q=1$. It yields $a_i \in L^{3-\delta}(\Omega),$ forall $i =1,\cdots,m-1$, where $\delta=\frac{9}{243}$. In other words there exists positive constant $C_{\delta_1}$, such that 
	
	\[
	\Vert a_i \Vert_{L^{\frac{234}{81}}(\Omega)}\leq C_{\delta_1}, \qquad \forall i=1,\cdots,m-1.
	\]
	 Next we will apply Gagliardo-Nirenberg interpolation. We obtain
	 
	 \[
	 \Vert a_i \Vert_{L^{3.9}(\Omega)} \leq C_{GN}\Vert \nabla a_i \Vert_{L^2(\Omega)}^{\frac{1}{2}} \Vert a_i \Vert_{L^{\frac{234}{81}}(\Omega)}^{\frac{1}{2}}+C_{GN}\Vert a_i \Vert_{L^{\frac{234}{81}}(\Omega)}, \qquad \forall i=1,\cdots,m-1,
	 \]
	 
	where the positive constant $C_{GN}$ depends only on domain and
	\[
	\frac{1}{3.9}=\bigg(\frac{1}{2}-\frac{1}{3}\bigg)\frac{1}{2}+\frac{81}{234}\bigg(1-\frac{1}{2}\bigg).
	\]

	Integrating with respect to time variable and applying H\"older inequality on right hand side, we obtain:
	
	\[
	\int_{0}^{t}\Vert a_i \Vert_{L^{3.9}(\Omega)}^{3.9} \leq 2^{2.9}C_{GN}^{3.9}C_{\delta_1}^{1.95}(1+T)\int_{0}^{T} \Vert \nabla a_i \Vert_{L^2(\Omega)}^2+2^{2.9}C_{GN}^{3.9}TC_{\delta_1}^{3.9}, \qquad \forall i=1,\cdots,m-1.
	\]
	
	Relations \eqref{L^2 gradient estimate} and relation \eqref{towards gradient estimate} yield 
	
	\begin{flalign*}
		\int_{0}^{t}\Vert a_i \Vert_{L^{3.9}(\Omega)}^{3.9} \leq 2^{2.9}C_{GN}^{3.9}(C_{\delta_1}^{1.95}+C_{\delta_1}^{3.9})(1+T) \bigg(\sum\limits_{i=0}^{m-1}\Vert a_{i,0} \Vert_{L^2(\Omega)}+C_{\Omega}(1+T)\bigg)
		 =C_2, \ \  \ \forall i=1,\cdots,m-1.
	\end{flalign*}
	
	Consider the o.d.e.,
	\[
	\partial_t a_m \leq \displaystyle\prod_{j=1}^{m-1} a_j^{\alpha_j}
	\] 
	Applying Minkwoski's inequality and H\"older inequality on the integral relation corresponding to above differential equation, we obtain
	
	\begin{flalign*}
		\Vert a_m \Vert_{L^{\frac{11.7}{10}}(\Omega)} \leq \Vert a_{m,0} \Vert_{L^{\frac{3.9}{3}}(\Omega)}+  \sum\limits_{i=1}^{m-1}\int_{0}^{t}\Vert a_i \Vert_{L^{3.9}(\Omega)}^{\frac{10}{3}} 
		 \leq \Vert a_{m,0} \Vert_{L^{\frac{3.9}{3}}(\Omega)}+mC^{\frac{100}{117}}_2T^{\frac{17}{117}}=C_3.
	\end{flalign*}
	
	From relation \eqref{Lp-Lq estimate} again, we take $q=\frac{11.7}{10},p=\frac{116}{22}$. It yields $a_i \in L^{\frac{116}{22}}(\Omega)$, for all $i=1,\cdots,m-1$. Let's assume
	there exists constant $C_{\delta_2}>0$, such that $\Vert a_i \Vert_{L^{\frac{116}{22}}(\Omega)}\leq C_{\delta_2}$.
	
	Again applying Minkwoski's integral inequality and H\"older inequality on the integral relation corresponding to $\partial_t a_m \leq \displaystyle\prod_{j=1}^{m-1} a_j^{\alpha_j}$, we obtain
	
	\begin{flalign*}
		\Vert a_m \Vert_{L^{\frac{348}{220}}(\Omega)} \leq \Vert a_{m,0} \Vert_{L^{\frac{348}{220}}(\Omega)}+ \sum\limits_{i=1}^{m-1}\int_{0}^{t}\Vert a_i \Vert_{L^{\frac{116}{22}}(\Omega)}^{\frac{10}{3}}
		 \leq \Vert a_{m,0} \Vert_{L^{\frac{348}{220}}(\Omega)}+mTC_{\delta_2}^{\frac{10}{3}}=C_4.
	\end{flalign*}
	
	As $\frac{348}{220}>\frac{3}{2}$, in relation \eqref{Lp-Lq estimate} we can simply put $q=\frac{348}{220}$ and $p=+\infty$, which provides

	\[
	\Vert a_i \Vert_{L^{\infty}(\Omega_T)} < +\infty, \qquad \forall i=1,\cdots,m-1.
	\]
	
	This further implies $a_m \in L^{\infty}(\Omega_T)$. Hence smooth unique global in time solution exists.

	\bibliography{ref3.bib}
	\bibliographystyle{abbrv}
	
	\appendix
	\section{some useful results}
	
	\begin{Thm}[$L^1$ Integrability estimate]\label{2nd order estimate}
		Let $d_1,\cdots,d_k\geq 0$ and $a_1(t,x),\cdots,a_k(t,x)\geq 0$, satisfy the following equation 
		\begin{flalign*}
			\partial_t \left(\sum\limits_{i=1}^{k} a_i\right)- \Delta \left(\sum\limits_{i=1}^{k} d_ia_i\right)\leq & 0 \qquad \qquad \mbox{\ in \ } \Omega_T, \ \forall i=1,\cdots,k\\
			n.\nabla_{x} a_i =& 0 \qquad \qquad \mbox{\ on \ } (0,T)\times \partial \Omega\\
			a_i(0,\cdot)(\geq 0) \in & L^2(\Omega) \qquad \mbox{\ in \ } \Omega.
		\end{flalign*}
		Then there exists a constant $\nu\geq 0$, depends on the initial data,domain and the dimension, such that 
		\[
		\int_{(0,T)\times \Omega} (\sum_{1}^{k} a_i)(\sum_{1}^{k} d_ia_i) \leq \nu.
		\]
	\end{Thm}
	Proof can be found in \cite{DFMV07}.
	
	\vspace{.2cm}
	
	\begin{Thm}[$L^p$ Duality estimate]\label{duality estimate}
		If  $\phi_1,\phi_2 \geq 0$, for some  $\Tilde{d_1},\Tilde{d_2}>0$, satisfy the following relation
		
		\begin{flalign*}
			\partial_t (\phi_1+\phi_2)-\Delta (\Tilde{d_1}\phi_1+\Tilde{d_1}\phi_2) \leq & 0 \qquad \ \ \ \Omega_T \\
			n.\nabla_{x}\phi_1=n.\nabla_{x}\phi_2= & 0 \qquad (0,T)\times \Omega \\
			\phi_1(0,x),\phi_2(0,x)\in & L^p(\Omega).
		\end{flalign*}
		
		Then  $\phi_1 \in L^p(\Omega_T) \implies \phi_2 \in L^p(\Omega_T)$  and vice versa.
	\end{Thm}
	
	Proof can be found in \cite{CDF14}\cite{DF08}.
	
	\vspace{.2cm}
	
	\begin{Thm}[Integrability estimation] \label{estimation 1}
		Let $d>0$, and let $\theta \in \mathrm L^p(\Omega_{T})$ for some $1<p< +\infty$. Let $\psi$ be the solution to the following parabolic equation
		\begin{flalign}
			\partial_t \psi(t,x) -d\Delta \psi(t,x)& = \theta(t,x) \qquad \quad \mbox{\ in \ } \Omega_{T} \nonumber\\
			n.\nabla_x \psi(t,x) & =0 \qquad \ \ \qquad \  \mbox{\ on \ } (0,T) \times \partial \Omega  \nonumber\\
			\psi(0,x) & =0 \qquad \qquad \ \ \  \ \mbox{\ in \ }\Omega. \nonumber
		\end{flalign}
		We have the following estimates
		\begin{equation}\label{eq:Ls-estimates}
			\begin{aligned}
				\Vert \psi \Vert_{L^s(\Omega _{T})} & \leq C_{IE} \Vert \theta \Vert_{L^p(\Omega_{T})} \quad\forall s< \frac{(N+2)p}{N+2-2p} \, \mbox{ with }\, p<\frac{N+2}{2},
				\\
				\Vert \psi \Vert_{L^s(\Omega _{T})} & \leq  C_{IE} \Vert \theta \Vert_{L^p(\Omega_{T})} \quad\forall s< +\infty \, \mbox{ with }\, p=\frac{N+2}{2},
			\end{aligned}
		\end{equation}
		where the constant $C_{IE}=C_{IE}(T,\Omega,d,p,s)$ and
		\begin{align}\label{eq:Linfty-estimate}
			\Vert \psi \Vert_{L^{\infty}(\Omega _{T})} \leq  C_{IE}(T, \Omega,d,p) \Vert \theta \Vert_{L^p(\Omega_{T})} \qquad  \mbox{ with } p> \frac{N+2}{2}.
		\end{align}
	\end{Thm}
	Proof of the estimates \eqref{eq:Ls-estimates} can be found in \cite[Lemma 3.3]{CDF14} and the estimate \eqref{eq:Linfty-estimate} was derived in \cite[Lemma 4.6]{Tan18}. We refer to the constant $C_{IE}$ as the integrability estimation constant.
\end{document}